\documentclass[12pt]{amsart}
\usepackage[cp1251]{inputenc}
\usepackage[english,russian]{babel}
\usepackage{amsmath}
\usepackage{amssymb}
\usepackage{amsfonts}
\title{On Some Transformations Associated to a Certain Cone}

\author{Vladimir B. Vasilyev}
\author{Denis A. Tokarev}

\date{}

\begin{document}
\renewcommand{\refname}{References}
\maketitle

{\bf Abstract.} A model elliptic pseudo-differential equation in $4$-faced cone is studied in Sobolev--Slobodetskii space. The Bochner kernel
for such a cone is evaluated and explicit formula for unique solution  to the considered equation is presented under certain restrictions on the symbol.
Boundary value problem with additional integral condition is considered and unique solvability to the boundary value problem is proved


\section{Introduction}

The theory of boundary value problems for elliptic equations has a long history, and there are a lot of publications in mathematical literature \cite{Ko,RS,S,NP,NSSS,Eg,SSS}.
New approach was suggested in \cite{V1}, it is related to a special factorization of an elliptic symbol at singular point of a manifold and reduces the general problem to invertibility conditions for so called "model" operators \cite{V4}.

The paper is related to some integral operators which are multidimensional analogues of the Hilbert transform \cite{G,M,MP,GK,K}. Such operators were used in \cite{V1} to construct solutions to an elliptic model pseudo-differential equation in a cone. One-dimensional cone is alone, it is half-axis, this fact was used in \cite{E} to construct solution in a half-space, but multidimensional space includes a lot of convex cones, and each such cone generates the Bochner kernel \cite{BM,Vl1,Vl2}. For convex cone $C\subset\mathbb R^m$ non including a whole straight line this kernel has the form
\begin{equation}\label{1}
	B_m(z)=\int\limits_Ce^{ix\cdot z}dx,
\end{equation}
where $z-(z_1,z_2,z_3),\,z_j-\xi_j+i\tau_j,\,j=1,2.3$.

The paper is devoted to a concrete cone in $\mathbb R^3$, evaluation of its Bochner kernel, and construction of a solution to model elliptic pseudo-differential equation and a
certain boundary value problem with an integral condition.

\section{$4$-faced cone in $\mathbb R^3$}

Here we consider the cone which is an intersection of five half-spaces:
\[
x_3>ax_1+bx_2,~~~x_3>cx_1+bx_2,~~~x_3>cx_1+dx_2,
\]
\[
x_3>ax_1+dx_2,~~~x_3>0,~~~a,b,c,d>0.
\]
The equation for such a conical surface is given by the formula
\[
x_3=\varphi(x_1,x_2)>0,
\]
and
$$
\varphi(x_1,x_2)=\begin{cases}
	ax_1+bx_2,~~~x_1>0,x_2>0,\\
	-cx_1+bx_2,~~~x_1<0,x_2>0,\\
	-cx_1-dx_20,~~~x_1<0,x_2<0,\\
	ax_1-dx_2,~~~x_1>0,x_2<0.
\end{cases}
$$

Thus, we consider the cone $C_3$ which is given by the inequality $x_3>\varphi(x_1,x_2)$. We will calculate the Bochner kernel for the cone $C_3$.

{\bf Lemma 1.}
{\it
	The Bochner kernel for $C_3$ has the following form
	\[
	B_3(z)=-\frac{i(a+c)(b+d)z_3}{(z_1+az_3)(z_2+bz_3)(z_1-cz_3)(z_2-dz_3)}.
	\]
}
{\bf Proof.}
	According to \eqref{1} we have
	\[
	B_3((\xi_1,\xi_2,z_3)=\int\limits_{C_3}e^{x_1\xi_1+x_2\xi_2+x_3z_3}dx=
	\]
	\[
	=\int\limits_0^{+\infty}e^{ix_1\xi_1}\left(\int\limits_0^{+\infty}e^{ix_2\xi_2}\left(\int\limits_{ax_1+bx_2}^{+\infty}e^{ix_3z_3}dx_3\right)dx_2\right)dx_1+
	\]
	\[
	+\int\limits_{-\infty}^{0}e^{ix_1\xi_1}\left(\int\limits_0^{+\infty}e^{ix_2\xi_2}\left(\int\limits_{cx_1+bx_2}^{+\infty}e^{ix_3z_3}dx_3\right)dx_2\right)dx_1+
	\]
	\[
	+\int\limits_{-\infty}^{0}e^{ix_1\xi_1}\left(\int\limits_{-\infty}^{0}e^{ix_2\xi_2}\left(\int\limits_{cx_1+dx_2}^{+\infty}e^{ix_3z_3}dx_3\right)dx_2\right)dx_1+
	\]
	\[
	+\int\limits_{0}^{+\infty}e^{ix_1\xi_1}\left(\int\limits_{-\infty}^{0}e^{ix_2\xi_2}\left(\int\limits_{bx_1+cx_2}^{+\infty}e^{ix_3z_3}dx_3\right)dx_2\right)dx_1=I_1+I_2+I_3+I_4.
	\]
	
	Let us evaluate the interior integral for $I_1$. We have
	\[
	\int\limits_{ax_1+bx_2}^{+\infty}e^{ix_3z_3}dx_3=-\frac{1}{iz_3}e^{iz_3(ax_1+bx_2)}=\frac{i}{z_3}e^{iaz_3x_1}e^{ibz_3x_2}.
	\]
	Further,
	\[
	\int\limits_0^{+\infty}e^{ix_2\xi_2}e^{ibz_3x_2}dx_2=\int\limits_0^{+\infty}e^{ix_2(\xi_2+bz_3)}dx_2=\frac{i}{\xi_2+bz_3}.
	\]
	And finally, the exterior integral
	\[
	\int\limits_{0}^{+\infty}e^{ix_1\xi_1}e^{iaz_3x_1}dx_1=\int\limits_{0}^{+\infty}e^{ix_1(\xi_1+az_3)}dx_1=\frac{i}{\xi_1+az_3},
	\]
	so that collection all evaluations we have
	\[
	I_1=-\frac{i}{z_3(\xi_1+az_3)(\xi_2+bz_3)}.
	\]
	
	For the integral $I_2-I_4$ all evaluations are almost the same and we will write shortly.
	\[
	\int\limits_{-cx_1+bx_2}^{+\infty}e^{ix_3z_3}dx_3=\frac{i}{z_3}e^{-icz_3x_1}e^{ibz_3x_2},
	\]
	the middle integral is the same,
	\[
	\int\limits_{-\infty}^{0}e^{ix_1\xi_1}e^{-icz_3x_1}dx_1=-\frac{i}{\xi_1-cz_3},
	\]
	so that
	\[
	I_2=\frac{i}{z_3(\xi_1-cz_3)(\xi_2+bz_3)},
	\]
	and
	\[
	I_1+I_2=-\frac{i}{z_3(\xi_1+az_3)(\xi_2+bz_3)}+\frac{i}{z_3(\xi_1-cz_3)(\xi_2+bz_3)}=
	\]
\[
=\frac{i(a+c)}{(\xi_1+az_3)(\xi_2+bz_3)(\xi_1-cz_3)}.
\]
	
	We continue with the second pair of summands. For $I_3$,
	\[
	\int\limits_{-cx_1-d-x2}^{+\infty}e^{ix_3z_3}dx_3=-\frac{1}{iz_3}e^{iz_3(-cx_1-dx_2)}=\frac{i}{z_3}e^{-icz_3x_1}e^{-idz_3x_2},
	\]
	the next integral on $x_2$
	\[
	\int\limits_{-\infty}^0e^{ix_2\xi_2}e^{-idz_3x_2}dx_2=\int\limits_{-\infty}^0e^{ix_2(\xi_2-dz_3)}dx_2=-\frac{i}{\xi_2-dz_3},
	\]
	and the third one on $x_1$
	\[
	\int\limits_{-\infty}^{0}e^{ix_1\xi_1}e^{-icz_3x_1}dx_1=\int\limits_{-\infty}^{0}e^{ix_1(\xi_1-cz_3)}dx_1=-\frac{i}{\xi_1-cz_3}.
	\]
	Hence,
	\[
	I_3=-\frac{i}{z_3(\xi_1-cz_3)(\xi_2-dz_3)}.
	\]
	
	Last, for $I_4$ we have
	\[
	\int\limits_{ax_1-d-x2}^{+\infty}e^{ix_3z_3}dx_3=-\frac{1}{iz_3}e^{iz_3(ax_1-dx_2)}=\frac{i}{z_3}e^{iaz_3x_1}e^{-idz_3x_2},
	\]
	the integral on $x_2$ as above for $I_3$,
	and the integral on $x_1$ is the following
	\[
	\int\limits_{0}^{+\infty}e^{ix_1\xi_1}e^{iaz_3x_1}dx_1=\int\limits_{0}^{+\infty}e^{ix_1(\xi_1+az_3)}dx_1=\frac{i}{\xi_1+az_3}.
	\]
	Thus,
	\[
	I_4=\frac{i}{z_3(\xi_1+az_3)(\xi_2-dz_3)}.
	\]
	So, we have
	\[
	I_3+I_4=-\frac{i}{z_3(\xi_1-cz_3)(\xi_2-dz_3)}+\frac{i}{z_3(\xi_1+az_3)(\xi_2-dz_3)}=
	\]
\[
=-\frac{i(a+c)}{(\xi_1+az_3)(\xi_2-dz_3)(\xi_1-cz_3)},
\]
therefore, collecting all integrals we conclude that
	\[
	B_3(\xi_1,\xi_2,z_3)=-\frac{i(a+c)(b+d)z_3}{(\xi_1+az_3)(\xi_1-cz_3)(\xi_2+bz_3)(\xi_2-d\xi_3)},
	\]
	and the proof is completed.$\blacksquare$

{\bf Remark 1.}
{\it	If $a=c, b=d$ then this case was considered in \cite{V5,V3}.}

\section{Pseudo-differential equations and wave factorization}

Let $C\subset\mathbb R^m$ be a sharp convex cone.
Principal equation under consideration is the following one
\begin{equation}\label{3}
	(Au)(x)=v(x),~~~x\in C,
\end{equation}
where $A$ is a model pseudo-differential operator with the symbol $A(\xi)$
satisfying the condition
\[
c_1(1+|\xi|)^{\alpha}\leq|A(\xi)|\leq c_2(1+|\xi|)^{\alpha},~~~c_1,c_2>0.
\]

Solution to the equation \eqref{3} is sought in the space
$H^s(C)$ is a subspace of $H^s(\mathbb R^m)$ with induced norm
\[
||u||s=\left(\int\limits_{\mathbb R^m}\tilde u(\xi)|^2(1+|\xi(^{2s}d\xi\right)^{1/2}
\]
and $supp~u\subset\overline{C}$, $\tilde u$ is the Fourier transform of $u$,
\[
\tilde u(\xi)=\int\limits_{\mathbb R^m}e^{ix\cdot\xi}u(x)dx.
\]
The operator $A$ is a bounded linear operator $A: H^s(\mathbb R^m)\rightarrow H^{s-\alpha}(\mathbb R^m)$, and the right hand side is taken from the space $\mathring{H}^{s-\alpha}(C)$; it consists of tempered distributions $S'(C)$ which admit continuation onto $H^{s-\alpha}(\mathbb R^m)$ with finite norm
\[
||v||^+_{s-\alpha}=\inf||\ell v||_{s-\alpha},
\]
where infimum is taken over all continuations $\ell$. 

For studying solvability of the equation \eqref{3} we use a special representation of an elliptic symbol.
Let us denote $T(C)=\mathbb R^m+iC$ radial tube domain over the cone $C$ and $\stackrel{*}{C}$ conjugate cone with respect to $C$ \cite{BM,Vl1,Vl2}.

{\bf Definition 1.}
{\it
	Wave factorization of elliptic symbol with respecr to the cone $C$ is called its representation in the form
	$$
	A(\xi)=A_{\neq}(\xi)A_=(\xi),
	$$
	where the factors $A_{\neq}(\xi),A_=(\xi)$ must satisfy the following conditions:
	
	1) $A_{\neq}(\xi),A_=(\xi)$ are defined for all admissible values $\xi\in{\mathbb R}^m$, without may be, the points $\{\xi\in{\mathbb R}^m:|\xi'|^2=a^2\xi^2_m\}$;
	
	2) $A_{\neq}(\xi),A_=(\xi)$ admit an analytical continuation into radial tube domains\\ $T(\stackrel{*} {C}),T(-\stackrel{*} {C})$ respectively with estimates
	$$
	|A_{\neq}^{\pm 1}(\xi+i\tau)|\leq c_1(1+|\xi|+|\tau|)^{\pm\ae},
	$$
	$$
	|A_{=}^{\pm 1}(\xi-i\tau)|\leq c_2(1+|\xi|+|\tau|)^{\pm(\alpha-\ae)},~\forall\tau\in\stackrel{*} {C}.
	$$
	
	The number $\ae\in{\mathbb R}$ is called index of wave factorization.}

Using such a factorization and results from \cite{V1} we can write explicit formula for the solution of equation \eqref{3} in $H^s{C_3}.$.

{\bf Theorem 1.}
{\it
	If the symbol $A(\xi)$ admits the wave factorization with respect to $C_3$ and $|\ae-s|<1/2$ then the equation \eqref{3} has unique solution for an arbitrary right hand side
	$v\in\mathring{H}^{s-\alpha}(C_3)$, and the solution is given by the formula
	\[
	\tilde u(\xi)=A^{-1}_{\neq}(\xi)\lim\limits_{\tau\to 0+}\int\limits_{\mathbb R^3}B_2(\xi'-\eta',\xi_3-\eta_3+i\tau)A^{-1}_=(\eta)\widetilde{(\ell v)}(\eta)d\eta,~~~\xi'=(\xi_1,\xi_2),
	\]
	where $\ell v$ is an arbitrary continuation of $v$ onto $H^{s-\alpha}(\mathbb R^3$.
	
	The a priori estimate
	\[
	||u||_s\leq~const~||v||^+_{s-\alpha}
	\]
	holds.}

\section{Transmutation operators}

According to construction of a solution to pseudo-differential equation \eqref{1} there is a certain transmutation operator
$V_{\varphi}$ \cite{V6}. The operator $V_{\varphi}$ is determined by geometry of a cone, and such an operator is distinct for different cones.
Here we will describe form of such an operator for two types of cones. Let us remind that the operator $T_{\varphi}: \mathbb R^m\rightarrow\mathbb R^m$
is the following
\begin{equation}\label{4}
	\begin{split}
		y_1=x_1\\
		y_2=x_2\\
		\cdots\\
		y_{m-1}=x_{m-1}\\
		y_m=x_m-\varphi(x_1,x_2,\dots,x_{m-1}
	\end{split}
\end{equation}

Let us introduce the following operators acting on different variables
\[
P_k=\frac{1}{2}(I+S_k),~~~~~Q_k=\frac{1}{2}(I-S_k),~~~k=1,2,
\]
where $I$ is identity operator,$S_k$ are one-dimensional singular integral operators of thr following type
\[
(S_1\tilde u)(\xi_1,\xi_2,\xi_3)=\frac{i}{\pi}v.p.\int\limits_{-\infty}^{+\infty}\frac{\tilde u(\eta,\xi_2,\xi_3)d\eta}{\xi_1-\eta},
\]
\[
(S_2\tilde u)(\xi_1,\xi_2,\xi_3)=\frac{i}{\pi}v.p.\int\limits_{-\infty}^{+\infty}\frac{\tilde u(\xi_1,\eta,\xi_3)d\eta}{\xi_2-\eta},
\]

Let us denote $\chi_{\pm}(x_k)=\{x\in\mathbb R^3: \pm x_k>0\},\,k=1,2$. It is well known \cite{E} that
\[
F_{x\to\xi}(\chi_{\pm}(x_k)u(x))=((I\pm S_k)\tilde u)(\xi)
\]
at least for $u\in S(\mathbb R^3)$, and we will use widely this property in evaluations below.

{\bf Lemma 2.} {\it
	Transmutation operator $V_{\varphi}=FT_{\varphi}F^{-1}$ acts as follows
	\[
	(V_{\varphi}\tilde u)(\xi)=
	(P_1P_2\tilde u)(\xi_1+a\xi_3,\xi_2+b\xi_3,\xi_3)+(Q_1P_2\tilde u)(\xi_1-c\xi_3,\xi_2+b\xi_3,\xi_3)+
	\]
	\[
	+(Q_1Q_2\tilde u)(\xi_1-c\xi_3,\xi_2-d\xi_3,\xi_3)+(P_1Q_2\tilde u)(\xi_1+a\xi_3,\xi_2-d\xi_3,\xi_3),
	\]
	at least for $u\in S{\mathbb R^3}$.
}

{\bf Proof.}
	So, we have
	\[
	(T_{\varphi}u)(x)=u(x_1,x_2,x_2-\varphi(x_1,x_2)),
	\]
	and then
	\[
	(FT_{\varphi}u)(\xi)=\int\limits_{\mathbb R^3}e^{ix\cdot\xi}u(x_1,x_2,x_3-\varphi(x_1,x_2))dx.
	\]
	We apply change of variables \eqref{4} and take into account that Jacobian
	$$
	\begin{vmatrix}
		1&0&0\\
		0&1&0\\
		-\frac{\partial\varphi}{\partial x_1}&-\frac{\partial\varphi}{\partial x_2}&1
	\end{vmatrix}
	$$
	equals to $1$ almost everywhere (it is piecewise smooth function). Then we have
	\[
	(FT_{\varphi}u)(\xi)=\int\limits_{\mathbb R^3}e^{iy_1\xi_1+iy_2\xi_2+i(y_3+\varphi(y_1,y_2))\xi_3}u(y_1,y_2,y_3)dy=
	\]
	\[
	=\int\limits_{\mathbb R^2}e^{iy_1\xi_1+iy_2\xi_2+i\varphi(y_1,y_2)\xi_3}\hat u(y_1,y_2,\xi_3)dy_1dy_2,
	\]
	where $\hat u$ is the Fourier transform on third variable.
	
	Further, we evaluate
	\[
	(FT_{\varphi}u)(\xi)=\int\limits_{\mathbb R^2}e^{iy_1\xi_1+iy_2\xi_2+i\varphi(y_1,y_2)\xi_3}\chi_+(y_1)\chi_+(y_2)\hat u(y_1,y_2,\xi_3)dy_1dy_2+
	\]
	\[
	+\int\limits_{\mathbb R^2}e^{iy_1\xi_1+iy_2\xi_2+i\varphi(y_1,y_2)\xi_3}\chi_-(y_1)\chi_+(y_2)\hat u(y_1,y_2,\xi_3)dy_1dy_2+
	\]
	\[
	+\int\limits_{\mathbb R^2}e^{iy_1\xi_1+iy_2\xi_2+i\varphi(y_1,y_2)\xi_3}\chi_-(y_1)\chi_-(y_2)\hat u(y_1,y_2,\xi_3)dy_1dy_2+
	\]
	\[
	+\int\limits_{\mathbb R^2}e^{iy_1\xi_1+iy_2\xi_2+i\varphi(y_1,y_2)\xi_3}\chi_+(y_1)\chi_-(y_2)\hat u(y_1,y_2,\xi_3)dy_1dy_2=
	\]
	\[
	=\int\limits_{\mathbb R^2}e^{iy_1\xi_1+iy_2\xi_2+i(ay_1+by_2)\xi_3}\chi_+(y_1)\chi_+(y_2)\hat u(y_1,y_2,\xi_3)dy_1dy_2+
	\]
	\[
	+\int\limits_{\mathbb R^2}e^{iy_1\xi_1+iy_2\xi_2+i(-cy_1+by_2)\xi_3}\chi_-(y_1)\chi_+(y_2)\hat u(y_1,y_2,\xi_3)dy_1dy_2+
	\]
	\[
	+\int\limits_{\mathbb R^2}e^{iy_1\xi_1+iy_2\xi_2+i(-cy_1-dy_2)\xi_3}\chi_-(y_1)\chi_-(y_2)\hat u(y_1,y_2,\xi_3)dy_1dy_2+
	\]
	\[
	+\int\limits_{\mathbb R^2}e^{iy_1\xi_1+iy_2\xi_2+i(ay_1-dy_2)\xi_3}\chi_+(y_1)\chi_-(y_2)\hat u(y_1,y_2,\xi_3)dy_1dy_2.
	\]
	Simple transformations lead to the following representation
	\[
	(FT_{\varphi}u)(\xi)=\int\limits_{\mathbb R^2}e^{iy_1(\xi_1+a\xi_3)+iy_2(\xi_2+b\xi_3)}\chi_+(y_1)\chi_+(y_2)\hat u(y_1,y_2,\xi_3)dy_1dy_2+
	\]
	\[
	+\int\limits_{\mathbb R^2}e^{iy_1(\xi_1-c\xi_3)+iy_2(\xi_2+b\xi_3)}\chi_-(y_1)\chi_+(y_2)\hat u(y_1,y_2,\xi_3)dy_1dy_2+
	\]
	\[
	+\int\limits_{\mathbb R^2}e^{iy_1(\xi_1-c\xi_3)+iy_2(\xi_2-d\xi_3)}\chi_-(y_1)\chi_-(y_2)\hat u(y_1,y_2,\xi_3)dy_1dy_2+
	\]
	\[
	+\int\limits_{\mathbb R^2}e^{iy_1(\xi_1+a\xi_3)+iy_2(\xi_2-d\xi_3)}\chi_+(y_1)\chi_-(y_2)\hat u(y_1,y_2,\xi_3)dy_1dy_2=
	\]
	\[
	=(P_1P_2\tilde u)(\xi_1+a\xi_3,\xi_2+b\xi_3,\xi_3)+(Q_1P_2\tilde u)(\xi_1-c\xi_3,\xi_2+b\xi_3,\xi_3)+
	\]
	\[
	+(Q_1Q_2\tilde u)(\xi_1-c\xi_3,\xi_2-d\xi_3,\xi_3)+(P_1Q_2\tilde u)(\xi_1+a\xi_3,\xi_2-d\xi_3,\xi_3),
	\]
	and Lemma 2 is proved.$\blacksquare$

{\bf Corollary.} {\it
	$$
	(V_{\varphi}\tilde u)(\xi_1,\xi_2,0)=\tilde u(\xi_1,\xi_2,0).
	$$
}

{\bf Remark 1.}
{\it
	If $\tilde u(\xi)=\tilde c(\xi'),\,\xi'=(\xi_1,\xi_2)$ then $(V_{\varphi}\tilde u(\xi_1,\xi_2,0)=\tilde c(\xi_1,\xi_2)$. Moreover, one can note that
	$V^{-1}_{\varphi}=V_{-\varphi}$.
}

\section{Boundary value problems}

We consider here the homogeneous equation
\begin{equation}\label{10}
	(Au)(x)=0,~~~x\in C_3
\end{equation}
with following integral condition
\begin{equation}\label{11}
	\int\limits_{-\infty}^{+\infty}u(x',x_3)=f(x_1,x_2).
\end{equation}

{\bf Theorem 2.}
{\it
	If the symbol $A(\xi)$ admits the wave factorization with respect to $C_3$ with the index $\ae$
	such that
	\begin{equation}\label{12}
		\ae-s=1+\varepsilon,\,|\varepsilon|<1/2
	\end{equation}
	then the boundary value problem \eqref{10},\eqref{11} has unique solution for arbitrary $f\in H^{s+1/2}(\mathbb R^2)$,
	and the Fourier transform of the solution can be represented in the form
	\[
	\tilde u(\xi)=(V^{-1}_{\varphi}\tilde c_0)(\xi_1,\xi_2,\xi_3),
	\]
	where $c_0\in H^{s-\ae+1/2}(\mathbb R^2)$ and
	\[
	\tilde c_0(\xi_1,\xi_2)=A(\xi_1,\xi_2,0)\tilde f(xi_1,\xi_2).
	\]
}

{\bf Proof.}
	In general, the proof is similar \cite{V5,V6} and we present main steps. First, we introduce the function $w(x)$ such that
	\[
	w(x)=\begin{cases}
		~~~~0,~~~~~~~~~~x\in C_3\\
		-(Au)(x),~~~x\in\mathbb R^3\setminus\overline{C_3},
	\end{cases}
	\]
	and rewrite the equation \eqref{10} in the form
	\[
	(Au)(x)+w(x)=0,~~~x\in\mathbb R^3.
	\]
	Then we apply the Fourier transform and the wave factorization and obtain the equality
	\[
	A_{\neq}(\xi)\tilde u(\xi)=-A^{-1}_=(\xi)\tilde w(\xi).
	\]
	Let us denote $\widetilde H(C)$ the Fourier image of the space $H(C)$. According to \cite{V1} we conclude that
	$A_{\neq}(\xi)\tilde u(\xi)\in\widetilde H^{s-\ae}(C_3),\,A^{-1}_=(\xi)\tilde w(\xi)\in\widetilde H^{s-\ae}(\mathbb R^3\setminus\overline{C_3})$.
	It means that the function $F^{-1}_{\xi\to 0}(A_{\neq}(\xi)\tilde u(\xi))$ is supported on $\partial C_3$ and after applying $t_{\varphi}$ it should be supported
	on the plane $x_3=0$. Such a distribution should be a span of summands \cite{E,Vl2} $c_k(x_1,x_2)\delta^{(k)}(x_3),\, k=0,1,\dots,n$, and each summand should be inside of $H^{s-\ae}(\mathbb R^3)$.
	According to condition \eqref{12} we have one summand only for $n=0$ so that
	\[
	T_{\varphi}F^{-1}_{\xi\to 0}(A_{\neq}(\xi)\tilde u(\xi))=c_0(x_1,x_2)\delta(x_3),
	\]
	where $\delta$ is Dirac mass-function. Estimates show that $c_0(x_1,x_2)\delta(x_3)\in H^{s-\ae}(\mathbb R^3)$ iff $c_0\in H^{s-\ae+1/2}(\mathbb R^2)$.
	
	After applying the Fourier transform we have
	\[
	V_{\varphi}(A_{\neq}(\xi)\tilde u(\xi)=\tilde c_0(\xi_1,\xi_2),
	\]
	in other words,
	\[
	\tilde u(\xi)=A^{-1}_{\neq})\xi)(V^{-1}_{\varphi}\tilde c_0)(\xi).
	\]
	
	It is left to determine arbitrary function $c_0$. We use the condition \eqref{11} which looks as follows
	\[
	\tilde u(\xi_1,\xi_2,0)=\tilde f(\xi_1,\xi_2).
	\]
	Then we have
	\[
	A_{\neq}(\xi_1,\xi_2,0)\tilde u(\xi_1,\xi_2,0)=\tilde c_0(\xi_1,\xi_2),
	\]
	or
	\[
	\tilde c_0(\xi_1,\xi_2)=A(\xi_1,\xi_2,0)\tilde f(xi_1,\xi_2),
	\]
	and Theorem 2 is proved.$\blacksquare$

{\bf Remark 2.}
{\it
	We have considered homogeneous equation \eqref{10}, but it was done for a simplicity. One can consider a non-homogeneous equation
	with the same methods also, but formulas obtained will be more large.
}

\section*{Conclusion}

This paper is related to special cone $C_3$ for which the Bochner kernel and corresponding transmutation operator were calculated. Moreover, the
boundary value problem from Sec. 5 is very specific and it is related to form of a general solution. May be it is possible to use other boundary condition to obtain explicit form for a solution but it is not clear to this moment.


\vspace{5mm}


Vladimir B. Vasilyev\\
Pobedy 85\\
Belgorod State National Research  University\\
Belgorod 308015\\
Russia\\
e-mail: vbv57@inbox.ru

\vspace{3mm}

Denis A. Tokarev\\
Pobedy 85\\
Belgorod State National Research  University\\
Belgorod 308015\\
Russia\\
e-mail: 1469493@bsuedu.ru\\

\end{document}